---

# Methods for arriving at numerical solutions for equations of the type (k+3) & (k+5) bi-quadratic's equal to a bi-quadratic (For different values of k )

Seiji Tomita[1] and Oliver Couto[2]

1) *Tokyo software company*, Tokyo, Japan

Email: fermat@m15.alpha-net.ne.jp

2) University of Waterloo, Ontario, Canada

Email: samson@celebrating-mathematics.com

Consider the below mentioned equations:

$$a^4 + b^4 + c^4 + d^4 + e^4 + k*(f)^4 = (g)^4 \;----(A)$$

$$For\ k = 1,2,3,4,5,6,7,8,9$$

$$a^4 + b^4 + c^4 + k*(d)^4 = (e)^4 \;-------(B)$$

$$For\ k = 2,3,7,8\ \&\ 9$$

**Abstract**: Different authors have done analysis regarding sums of powers (Ref. no. 1,2 & 3), but systematic approach for solving Diophantine equations having sums of many bi-quadratics equal to a quartic has not been done before. In this paper we give methods for finding numerical solutions to equation (A) given above in section one. Next in section two, we give methods for finding numerical solutions for equation (B) given above. As is known that finding parametric solutions to biquadratic equations is not easy by conventional method. So the authors have found numerical solutions to equation (A) & (B) using elliptic curve theory.

---

**Summary:**

Previously research has been done by others, for Sum of bi-quadratics equal to a bi-quadratic.

R. Carmichael gave the below mentioned identity for (k+2) bi-quadratics where k=4,

$$(a^4 - 2b^4)^4 + (2a^3b)^4 + 4*(2ab^3)^4 = (a^4 + 2b^4)^4$$

Also equation ( $a + k * b^4$=$c^4 + k * d^4$), for k=4 , has been investigated by Ajai Choudhry in his paper, but our paper is different, because of our paper's requirement, that the right hand side of our above equations (A) &(B) need's to be equal to a quartic. But to-date the authors have not found any publications which deal systematically, with the subject of (k + n) bi-quadratics, where (n=3, 5). We have, numerical solutions where k>9, but the authors are going to deal with solutions for k=1,2,3,4,5,6,7,8,9 for (k+5) bi-quadratics and for k=2,3,7,8,9 for (k+3) bi-quadratics . Mention of k=4 & 6 is not made in section two because k=4 & 6 has elliptic curve, with rank zero, hence has no rational integer solutions. K=5 is omitted in the latter because it has alluded solution by our method. Also k=1 is omitted because a solution has been given by Jacobi & Madden in their paper on the equation,

$$(a, b, c, d)^4 = (a + b + c + d)^4$$

**Section one:**

Equation $A^4 + B^4 + C^4 + D^4 + E^4 + k * F^4 = G^4$ has infinitely many integer solutions for k=1,2,3,4,5,6,7,8,9.

Proof.

$$A^4 + B^4 + C^4 + D^4 + E^4 + k * F^4 = G^4 \text{.............................. (1)}$$

$$Let\ A = c * x^2 + d, B = e * x + f, C = s * r, D = t * r, E = u * r, F = r,$$

$$G = a * x^2 + b \ldots\ldots\ldots\ldots\ldots\ldots\ldots\ldots\ldots\ldots. (2)$$

$$Consider\ the\ equation,$$

$$r^4 * (s^4 + t^4 + u^4 + k) = (G^4 - A^4 - B^4) \ldots\ldots\ldots\ldots\ldots\ldots. (3)$$

$$G^4 - A^4 - B^4 = (\ ax^2 + b)^4 - (\ cx^2 + d)^4 - (\ ex + f)^4$$

= $(\ 8 + 8x + 32x^3)$^2 $for$ $[\mathrm{a, b, c, d, e, f}] = [4, 3, 4, -1, 4, -2]$, So Right hand side is a square and we need to make the left hand side of equation (3) $[r^4 * (s^4 + t^4 + u^4 + k)]\ a\ square$, for suitable values of (s,t,u) & (k)

We have below equation (A)

$$A^4 + B^4 + C^4 + D^4 + E^4 + k * F^4 = G^4$$

**Equation for, k=1**:

$$A^4 + B^4 + C^4 + D^4 + E^4 + F^4 = G^4$$

Let [s,t,u]=[19, 17, 11] and [a,b,c,d,e,f]=[4, 3, 4, -1, 4, -2],

$$equation\ (3)\ above\ , becomes$$

$$-4(-239r^2 + 4 + 4x + 16x^3)(239r^2 + 4 + 4x + 16x^3) = 0.$$

$$Hence\ -239r^2 + 4 + 4x + 16x^3 = 0\ or$$

$$239r^2 + 4 + 4x + 16x^3 = 0.$$

$$Let\ us\ find\ the\ rational\ solution\ of$$

$$-239r^2 + 4 + 4x + 16x^3 = 0 \dots\dots\dots\dots\dots\dots\dots\dots\dots\dots..(4)$$

$$Substitute\ x = \frac{X}{956}\ and\ r = \frac{Y}{114242}\ to\ equation\ (4), then\ we\ obtain$$

$$Y^2 = X^3 + 228484X + 218430704 \dots\dots\dots\dots\dots\dots\dots\dots\dots\dots\dots\dots\dots.(5)$$

$$The\ elliptic\ curve\ Y^2 = X^3 + 228484X + 218430704\ has\ rank\ 1.$$

Refer to the elliptic curve tables mentioned in the reference section:

$$The\ point\ (X,Y)\ on\ the\ elliptical\ curve\ above = (580\ , 23368)\ , leads\ to\ below.$$

$$26979^4 + 24378^4 + 221996^4 + 198628^4 + 128524^4 + 11684^4 = 255463^4$$

-------

$$We\ have,\quad A^4 + B^4 + C^4 + D^4 + E^4 + k * F^4 = G^4$$

**Equation for, k=2**.

$$A^4 + B^4 + C^4 + D^4 + E^4 + 2 * F^4 = G^4$$

Let [s,t,u]=[7, 3, 2] and [a,b,c,d,e,f]=[4, 3, 4, -1, 4, -2],

Equation (3) becomes,

$$-4(-25r^2 + 4 + 4x + 16x^3) * (25r^2 + 4 + 4x + 16x^3) = 0.$$

In the same way with k=1, let us find the rational solution of

$$-25r^2 + 4 + 4x + 16x^3 = 0 \text{........ (6)}$$

Let us transform equation (6) to elliptic curve (7).

$$Y^2 = X^3 + 4X + 16 \text{........................................................................ (7)}$$

The elliptic curve eqn.(7) has rank 1.

The point (X, Y) = (1/4, -33/8) leads to below mentioned numerical solution.

General Note: The elliptical equations in section one & two are all having rank one. Since the points (X,Y) on the elliptical curves are rational points, by Nagell-Lutz theorem the elliptical curves will all have points of infinite order.

$$315^4 + 560^4 + 924^4 + 396^4 + 264^4 + 2 * 132^4 = 965^4$$

$$A^4 + B^4 + C^4 + D^4 + E^4 + 2 * F^4 = G^4 \text{ --------(7a)}$$

Additional method for ‘k’=2

Consider the above equation (7a) for (k+5) quartics,

For k=2, there is a numerical solution for the above namely,

$$(6,10,16,32,29)^4 + 2*(12)^4 = (37)^4$$

[since (6+10=16)]

Example of parametric solution for above for k=2, meaning (a+b=c) is given below:

$$(6n^2 — 20n - 16)^4 + (-16n^2 - 12n = 10)^4 + (-10n^2 - 32n - 6)^4 + 32^4 * (n^2 + n + 1)^4 + 29^4 * (n^2 + n + 1)^4 + 2 * (n^2 + n + 1)^4 = (37n^2 + 37n + 37)^4$$

Where ‘n’ is a parameter.

$$(a, b, c, d, e)^4 + 2 * (f)^4 = g^4$$

As is known (a+b=c) $implies\ ((a, b, (a+b))^4 = 2 * (a^2 + ab + b^2)^2$, hence equation (1) above can be parametrized when (a,b) is known.

$$Since\ (a, b, c)^4 = 2(a^2 + ab + b^2)^2 = [g^4 - d^4 - e^4 - k * (f)^4]$$

---

$$We\ have,\ A^4 + B^4 + C^4 + D^4 + E^4 + k * F^4 = G^4$$

**Equation for, k=3.**

$$A^4 + B^4 + C^4 + D^4 + E^4 + 3 * F^4 = G^4$$

Let [s,t,u]=[5, 4, 2] and [a,b,c,d,e,f]=[4, 3, 4, -1, 4, -2],

Equation (3) becomes to

$$-4(15r^2 + 4 + 4x + 16x^3) * (-15r^2 + 4 + 4x + 16x^3) = 0.$$

Let us find the rational solution of $(-15r^2 + 4 + 4x + 16x^3)$= 0 ................................. (8)

Let us transform equation (8) to elliptic curve (9).

$$Y^2 = X^3 + 900X + 54000 \quad \text{.................................................................. (9)}$$

The elliptic curve equation (9) has rank 1.

The point (X, Y) = (34 , -352) leads to below.

$$16^4 + 15^4 + 220^4 + 176^4 + 88^4 + 3 * 44^4 = 241^4$$

$$A^4 + B^4 + C^4 + D^4 + E^4 + k * F^4 = G^4$$

Equation for, k=4.

$$A^4 + B^4 + C^4 + D^4 + E^4 + 4 * F^4 = G^4$$

Let [s,t,u]=[70, 30, 20] and [a,b,c,d,e,f]=[4, 3, 4, -1, 4, -2],

Equation (3) becomes ,

$$-4(2499r^2 + 4 + 4x + 16x^3) * (-2499r^2 + 4 + 4x + 16x^3) = 0.$$

Let us find the rational solution of

$$-2499r^2 + 4 + 4x + 16x^3 = 0 \text{ ............................. } (10)$$

Let us transform equation (10) to elliptic curve eqn.(11).

$$Y^2 = X^3 + 10404X + 2122416 \text{ .................................................... } (11)$$

The elliptic curve (11) has rank 1.

The point (X, Y) = ( 474 , 10656) leads to below.

$$10416^4 + 3689^4 + 10360^4 + 4440^4 + 2960^4 + 4 * 148^4 = 12439^4$$

---

$$We\ have, \quad A^4 + B^4 + C^4 + D^4 + E^4 + k * F^4 = G^4$$

**Equation for, k=5.**

$$A^4 + B^4 + C^4 + D^4 + E^4 + 5 * F^4 = G^4$$

Let [s,t,u]=[11, 7, 5] and [a,b,c,d,e,f]=[4, 1, 4, -1, 4, 0],

Equation (3) becomes,

$8 * (47r^2 + 2x - 8x^3) * (47r^2 - 2x + 8x^3) =0.$

Let us find the rational solution of $47r^2 + 2x - 8x^3$=0................................ (12)

Let us transform equation (12) to elliptic curve eqn.(13).

$Y^2 = X^3 - 2209X$ ..................................................................... (13)

The elliptic curve (13) has rank 1.

The point (X, Y) = (684407232/2289169, 17682275119320/3463512697)

leads to below numerical solution.

$$206807355454175^4 + 66669098675328^4 + 133221414581640^4 + 84777263824680^4 + 60555188446200^4 + (5) * 12111037689240^4 = 217287944875297^4$$

$$Next\ we\ have, \quad A^4 + B^4 + C^4 + D^4 + E^4 + k * F^4 = G^4$$

**Equation for, k=6.**

$$A^4 + B^4 + C^4 + D^4 + E^4 + 6 * F^4 = G^4$$

Let [s,t,u]=[37, 31, 11] and [a,b,c,d,e,f]=[4, 3, 4, -1, 4, -2],

Equation (3) becomes,

$$-(-1677r^2 + 8 + 8x + 32x^3) * (1677r^2 + 8 + 8x + 32x^3) = 0.$$

Let us find the rational solution of $-1677r^2 + 8 + 8x + 32x^3 = 0$ ............................ (14)

Let us transform equation (14) to elliptic curve (15).

$Y^2 = X^3 + 44997264X + 603683293824$ ..................................................... (15)

The elliptic curve (15) has rank 1.

The point (X, Y) = ( 1720 , 828352) leads to below.

$$1421^4 + 2262^4 + 4144^4 + 3472^4 + 1232^4 + 6 * 112^4 = 4663^4$$

$$A^4+B^4+C^4+D^4+E^4+k*F^4=G^4$$

Equation for, k=7.

$$A^4+B^4+C^4+D^4+E^4+7*F^4=G^4$$

Let [s,t,u]=[5, 3, 2] and [a,b,c,d,e,f]=[4, 3, 4, -1, 4, -2],

Equation (3) becomes,

$-(27r^2+8+8x+32x^3)(-27r^2+8+8x+32x^3)$=0.

Let us find the rational solution of

$-27r^2+8+8x+32x^3$= 0 ................................ (16)

Let us transform equation (16) to elliptic curve (17).

$Y^2 \ = \ X^3+144X+3456$ ................................................................. (17)

The elliptic curve (17) has rank 1.

The point (X, Y) = ( 4 , -64) leads to below.

$$6^4+\ 9^4+\ 20^4+\ 12^4+\ 8^4+\ 7*\ 4^4=\ 21^4$$

---

$$We\ have, \quad A^4+B^4+C^4+D^4+E^4+k*F^4=G^4$$

**Equation for, k=8.**

$$A^4+B^4+C^4+D^4+E^4+8*F^4=G^4$$

Let [s,t,u]=[4, 3, 2] and [a,b,c,d,e,f]=[4, 3, 4, -1, 4, -2],

Equation (3) becomes,

$-(-19r^2+8+8x+32x^3)*(19r^2+8+8x+32x^3)$ =0.

Let us find the rational solution of

$-19r^2+8+8x+32x^3$= 0 ................................ (18)

Let us transform equation (18) to elliptic curve (19).

$Y^2 = \ X^3+5776X+877952$ ............................................................ (19)

The elliptic curve eqn.(19) has rank 1.

The point (X, Y) = (-16316/225, -941248/3375) leads to below.

$$409346^4+17856675^4+3529680^4+2647260^4+1764840^4+\ 8*882420^4 = 17866279^4$$

---

$$We\ have, \quad A^4+B^4+C^4+D^4+E^4+k*F^4=G^4$$

**Equation for, k=9.**

$$A^4 + B^4 + C^4 + D^4 + E^4 + 9 * F^4 = G^4$$

Let [s,t,u]=[12, 10, 6] and [a,b,c,d,e,f]=[4, 3, 4, -1, 4, 2],

Equation (3) becomes ,

$$-(179r^2 - 8 + 8x + 32x^3) * (-179r^2 - 8 + 8x + 32x^3) = 0.$$

Let us find the rational solution of -179r^2-8+8x+32x^3=0 ............................... (20)

Let us transform equation (20) to elliptic curve eqn.(21).

$Y^2 = X^3 + 512656X - 734123392$ ..................................................... (21)

This elliptic curve eqn. (21) has rank 1.

The point (X, Y) = (5696705292406351213336/10878607024914721,

13598002320735074871580564215680/1134644815597146377458481)

leads to below numerical solution. ------------------------------------------- (22)

$$632907528785561577532579698212415075^4$$
$$+175473636600521434023931273346458144^4$$
$$+1327935398893889305717227119370758400^4$$
$$+1106612832411574421431022599475632000^4$$
$$+663967699446944652858613559685379200^4$$
$$+ 9 * 110661283241157442143102259947563200^4$$
$$= 633380905148771673201251847502446439^4$$

---

<u>**Section Two: (k+3) quartic's**</u>

Equation (B) is given below,

$$We\ have,\quad A^4 + B^4 + C^4 + k * D^4 = E^4 \ldots\ldots\ldots\ldots\ldots\ldots\ldots\ldots\ldots\ldots\ldots (1)$$

$$Let\ A = cx^2 + d, B = ex + f, C = sr, D = r, E = ax^2 + b \quad \ldots\ldots\ldots\ldots.. (2)$$

Hence we get the equation,

$$r^4 * (s^4 + k) = (E^4 - A^4 - B^4) \ldots\ldots\ldots\ldots\ldots\ldots\ldots\ldots.. (3)$$

$$G^4 - A^4 - B^4 = (a * x^2 + b)^4 - (c * x^2 + d)^4 - (e * x + f)^4$$

$= (8 + 8x + 32x^3)^2\ for$ $[a, b, c, d, e, f] = [4, 3, 4, -1, 4, -2]$, So Right hand side is a square and we need to make the left hand side of equation (3) $[r^4 * (s^4 + k)]\ a\ square$, for suitable values of (s) & (k)

Section Two:

Regarding the (k+3) quartic equation given below

**For k=2 ,**

$$(a,b,c)^4 + k * d^4 = e^4 \text{ ------ (B)}$$

For k=2, the Identity is given below,

$$(p^2+q^2)^4 = (p^2-q^2)^4 + (2pq)^4 + (s)^4 + 2 * (r)^4$$

The above after simplification has the condition:

$$3r^2 = 2pq(p^2-q^2) ----(C)$$

Solution is (p,q,r,s) = (2,1,2,4)

The above, equation (B) has the elliptical equation,

$$Y^2 = X * (X^2 - 36)$$

And it has rank one, hence equation (A) has infinite solutions for k=2.

Numerical solution is:

$$(4,3,4)^4 + 2 * (2)^4 = (5)^4$$

The point (X, Y) = (25/4, 35/8) leads to below numerical solution,

$$49^4 + 280^4 + 1200^4 + 2 * 140^4 = 1201^4$$

The author's have noted that there are parametric solutions for k= 14,46,49,63,94 where (k)<100. For (k+3) quartics. Those have the condition (a+b=c). Namely $(a,b,c)^4 + k * (d)^4 = e^4$

---

$$We\ have, \qquad A^4 + B^4 + C^4 + k * D^4 = E^4$$

**Equation for, k=3.**

$$A^4 + B^4 + C^4 + 3 * D^4 = E^4$$

Let [a,b,c,d,e,f,s]=[4, 3, 4, -1, 4, 2, 1/2],

Equation (3) becomes,

$$\frac{1}{16} * (7r^2 + 32 - 32x - 128x^3) * (7r^2 - 32 + 32x + 128x^3) = 0.$$

Let us find the rational solution of

$7r^2 + 32 - 32x - 128x^3 = 0$ ....................... (6)

Let us transform equation (6) to elliptic curve eqn.(7).

$Y^2 = X^3 + 784X - 43904$ ........................................................ (7)

The elliptic curve eqn. (7) has rank 1.

The point (X, Y) = (36, 176) leads to below.

$$8^4 + 56^4 + 11^4 + 3 * 22^4 = 57^4$$

$$Next\ we\ have, \quad A^4 + B^4 + C^4 + k * D^4 = E^4$$

**Equation for, k=7.**

$$A^4 + B^4 + C^4 + 7 * D^4 = E^4$$

When [a,b,c,d,e,f,s]=[4, 1, 4, -1, 4, 0, 47],

Equation (3) becomes,

$$8(781r^2 + 2x - 8x^3) * (781r^2 - 2x + 8x^3) = 0.$$

Let us find the rational solution of

$781r^2 + 2x - 8x^3 = 0$........................... (8)

Let us transform equation (8) to elliptic curve eqn.(9).

$Y^2 = X^3 - 609961X$ .......................................................... (9)

The elliptic curve (9) has rank 1.

The point (X, Y) = (-2876843001196439/4324112302500, -94873842643707990383059/8991775327433625000) leads to below.

$$5129496674953832213892839^4 + 3185606200725875569549500 0^4$$
$$+ 1520165120067767166801885 0^4 + 7 * 3234393872484610993195 50^4$$
$$= 3226639773430987079860716 1^4$$

---

$$We\ have, \quad A^4 + B^4 + C^4 + k * D^4 = E^4$$

**Equation for, k=8.**

$$A^4 + B^4 + C^4 + 8 * D^4 = E^4$$

Let [a,b,c,d,e,f,s] = [4, 3, 4, -1, 4, -2, 239/13],

Equation (3) becomes,

$-\frac{1}{28561} * (-57123r^2 + 1352 + 1352x + 5408x^3) * (57123r^2 + 1352 + 1352x + 5408x^3) = 0.$

Let us find the rational solution for,

$-57123r^2 + 1352 + 1352x + 5408x^3 = 0$.............(10)

Let us transform equation (10) to elliptic curve (11).

$Y^2 = X^3 + 18409008087184X + 1579703492932904584 96$ ...........(11)

The elliptic curve equation (11) has rank 1.

The point (X, Y) = (2088556756/1369, 697479284591232/50653) leads to below.

$$136268507232^4 + 201049446673^4 + 483363968776^4 + 8 * 26291763992^4 = 487694040337^4$$

$$A^4 + B^4 + C^4 + k * D^4 = E^4$$

Equation for, k=9.

$$A^4 + B^4 + C^4 + 9 * D^4 = E^4$$

When [a,b,c,d,e,f,s]=[4, 3, 4, -1, 4, 2, 2],

Equation (3) becomes,

$-(-5r^2 - 8 + 8x + 32x^3) * (5r^2 - 8 + 8x + 32x^3) = 0.$

Let us find the rational solution of ,

$-5r^2 - 8 + 8x + 32x^3 = 0$......................... (12)

Let us transform equation (12) to elliptic curve eqn.(13).

$Y^2 = X^3 + 400X - 16000$......................................................... (13)

This elliptic curve eqn.(13) has rank 1.

The point (X, Y) = (164, -2112) leads to below.

$$414^4 + 115^4 + 264^4 + 9 * 132^4 = 439^4$$

Similarly we also have numerical solutions for (k+5) quartics, in which C=(A+B)

$$A^4 + B^4 + C^4 + D^4 + E^4 + k * F^4 = G^4$$

Hence,

$$A^4 + B^4 + (A + B)^4 + D^4 + E^4 + k * F^4 = G^4$$

$(a, b, c)^4 + k * d^4 = e^4$ ------ (A).

$$Regarding\ the\ above\ equation\ (A), for\ the\ case\ (c = a + b).$$

$$Consider\ the\ below\ mentioned\ identity\ for\ (k = 4) and\ with\ the$$

$$condition\ (c = a + b). Let\ k = (m)^2 = (2)^2 = 4$$

$$(2p^2 - 2q^2)^4 + (2q^2 + 4pq)^4 + (2p^2 + 4pq)^4 + (4) * (p^2 + pq + q^2)^4 = [6 * (p^2 + pq + q^2)^2]^2$$

The right hand side of the above equation is equal to 36 times a fourth power and hence cannot be made a fourth power.

Regarding the case $k = (m)^2$, m=1,2,3,4,5,6 has no known solutions with the condition (c=a+b). But there is a solution for m=7, meaning 'k' = $(m)^2 = (7)^2 = 49$

Consider the below equation for (k+5) quartics,

$(a,b,c,d,e)^4 + (k) * (f)^4 = (g)^4$ --------- (A)

For k=5, there is a numerical solution for the above namely,

$$(4,22,26,7,28)^4 + 5 * (14)^4 = (35)^4$$

[since (4+22=26)]

Example of parametric solution for k=5, for (a+b=c) is given below:

$$(-26n^2 — 44n + 4)^4 + (22n^2 - 8n - 26)^4 + (4n^2 + 52n + 22)^4 + 5 * (n^2 + n + 1)^4 = (35n^2 + 35n + 35)^4$$

Where ‘n’ is a parameter.

There are parametric solutions for k= 14,46,49,63,94 where (k)<100. For (k=3) quartics. These have the condition (a+b=c). Namely $(a,b,c)^4 + k * (d)^4 = e^4$ --------- (1)

As is known (a=b=c) $implies\ ((a,b,(a+b))^4 = 2 * (a^2 + ab + b^2)^2$, hence equation (1) above can be parametrized when (a,b) is known.

$$Since\ (a,b,c)^4 = 2(a^2 + ab + b^2)^2 = e^4 - k * (d)^4$$

Example of parametric solution for k=14 is given below:

Numerical solution for above is $(4,11,15)^4 + 14 * (1)^4 = 16^4,$ (a+b=c) implies in above (4+11=15)

-------------

See below for more numerical solutions regarding (k+3) & (k+5) quartics,

for k=1,2,3,4,5,6,7,8 & 9

Table (1) Numerical solutions to (k+3) equation given below,

$$a^4 + b^4 + c^4 + k * (d)^4 = (e)^4$$

For k=1,2,3,4,5,6,7,8 & 9

| k | a | b | c | $d^*$ | e |
|---|---|---|---|---|---|
| 1 | 30 | 120 | 272 | 315 | 353 |
| 2 | 49 | 280 | 1200 | 140 | 1201 |
| 3 | 2 | 4 | 7 | 6 | 9 |
| 4 | 34 | 10 | 5 | 14 | 35 |
| 5 | 69 | 40 | 40 | 94 | 143 |
| 6 | 455 | 280 | 142 | 170 | 483 |
| 7 | 4 | 4 | 1 | 2 | 5 |
| 8 | 3 | 2 | 2 | 22 | 37 |
| 9 | 15 | 14 | 6 | 34 | 59 |

(*) means $(d^4)$ needs multiplication by value of 'k'

Conclusion:

The authors have provide d method's for finding numerical solutions for (k+3) & (k+5) quartic equation for k = (1,2,3,4,5,6,7,8 & 9) . Except for the case of (k=4,5,6) for the (k+3) quartic equation it is an open question and others can attempt to find solution by a different method.

Table (2) Numerical solutions for (K+5) quartic equation given below:

For k=1,2,3,4,5,6,7,8 & 9

$$a^4 + b^4 + c^4 + d^4 + e^4 + k * (f)^4 = (g)^4$$

| k | a | b | c | d | e | $f^*$ | g |
|---|---|---|---|---|---|---|---|
| 1 | 6 | 8 | 18 | 31 | 32 | 34 | 43 |
| 2 | 2 | 6 | 8 | 13 | 20 | 4 | 21 |
| 3 | 4 | 5 | 6 | 8 | 10 | 8 | 13 |
| 4 | 10 | 12 | 14 | 15 | 20 | 2 | 23 |
| 5 | 3 | 4 | 6 | 8 | 14 | 6 | 15 |
| 6 | 1 | 8 | 12 | 14 | 16 | 4 | 19 |
| 7 | 2 | 10 | 18 | 19 | 24 | 28 | 47 |
| 8 | 4 | 5 | 8 | 10 | 18 | 6 | 19 |
| 9 | 8 | 18 | 27 | 42 | 48 | 10 | 55 |

(*) means $(f^4)$ needs to be multiplied by value of (K)

Conclusion:

The authors have provided method's for finding numerical solutions for (k+5) quartic equation for k = (1,2,3,4,5,6,7,8 & 9). Also the authors have provided solutions for (k+3) quartic equation for (k+3) equation for k=2,3,7,8&9, and for the case of (k=4,5,6) it is an open question and others can attempt to find solution by a different method.

---